\documentclass[letterpaper]{article}
\usepackage[letterpaper,margin=1in]{geometry}

\usepackage{amsmath,amsfonts,amssymb,amsthm}

\usepackage[colorlinks,citecolor=red]{hyperref}

\newtheorem{definition}{Definition}[section]
\newtheorem{theorem}{Theorem}[section]
\newtheorem{lemma}{Lemma}[section]

\newtheorem{proposition}{Proposition}[section]

\newcommand{\RR}{\mathbb{R}}
\newcommand{\NN}{\mathbb{N}}

\newcommand{\norm}[1]{\left|#1\right|}

\DeclareMathOperator*{\esssup}{ess\,sup}

\title{Non-Existence of Solutions for a Non-Gaussian Equation in Fractional Time with Osgood Type Nonlinearity}

\author{Soveny Solís\thanks{Departamento de Matem\'atica, Facultad de Ciencias Naturales y Matem\'aticas, Escuela Superior Polit\'ecnica del Litoral, Guayaquil, Ecuador \\ \texttt{ssolis@espol.edu.ec}} \and Vicente Vergara\thanks{Departamento de Matem\'atica, Facultad de Ciencias F\'isicas y Matem\'aticas, Universidad de Concepci\'on,  Concepci\'on, Chile \\ \texttt{vvergaraa@udec.cl}}}
\date{}

\begin{document}

\maketitle


\abstract{Osgood functions in the source term are used to produce results for non-existence of local solutions into the framework of non-Gaussian diffusion equations. The critical exponent for non-existence of local solutions is found to depend on the fractional derivative, the non-Gaussian diffusion and the non-linear term. The instantaneous blow-up phenomenon is studied by exploiting estimates of the fundamental solutions. Nevertheless, theory of super-solutions and fixed points are combined for showing existence of global solutions. In this case, the critical exponent for existence of global solutions depends only on the last two parameters above.
} 


\begin{center}
\end{center}




\section{Introduction}


Non-Gaussian processes have gained increasing attention in recent years, as they are better suited for modeling complex systems that exhibit long-range dependence or heavy-tailed distributions. Fractional time derivatives, on the other hand, provide a powerful tool for modeling non-local phenomena, making them ideal for modeling a wide range of phenomena in physics, finance, biology, and other fields, see e.g. \cite{ADK23}, \cite{AwMe20}, \cite{Bel21}, \cite{LW21} and \cite{SaTo19}.

Consider $\alpha \in (0,1)$ and $\beta \in (0,2)$. We investigate the semilinear evolution problem:

\begin{equation}\label{FE:0}
	\partial_t^{\alpha} (u - u_0) + \Psi_{\beta}(-i\nabla) u = f(u), \quad t > 0, \quad x \in \mathbb{R}^d.
\end{equation}

Here, $\partial_t^{\alpha}$ denotes the Riemann-Liouville fractional derivative of order $\alpha$, defined as 
\[
\partial_t^{\alpha} v = \frac{d}{dt}\int_0^t g_{1-\alpha}(t-s)v(s)ds,
\] 
with $g_{\rho}(t) = \frac{1}{\Gamma(\rho)}t^{\rho-1}$ and $\Gamma(\cdot)$ denoting the gamma function. The operator $\Psi_{\beta}(-i\nabla)$ is a pseudo-differential operator of order $\beta$, characterized by the symbol $\psi_{\beta}(\xi) = |\xi|^{\beta} \omega_{\nu}(\xi/|\xi|)$, where $\omega_{\nu}$ is a continuous function on the surface of the $(d-1)$-dimensional sphere $\mathbb{S}^{d-1}$, and $\nu$ represents a spectral measure on $\mathbb{S}^{d-1}$. Specifically,

\begin{equation}\label{omega:nu}
	\omega_{\nu}(\theta) = \int_{\mathbb{S}^{d-1}} |\theta \cdot \eta|^{\beta}\nu(d\eta), \quad \theta \in \mathbb{S}^{d-1}.
\end{equation}
For more details on the operator $\Psi_{\beta}(-i\nabla)$, refer to \cite{Kolo09}.

\medbreak

The evolution problem \eqref{FE:0} is considered under the following hypotheses for the initial data $u_0$, the spectral measure $\nu$, and the source term $f$:
\begin{itemize}
	\item[(H1)] For $1\leq q < \infty$, $u_0\in L_q(\RR^d)$ and $u_0\geq 0$.
	\item[(H2)] The spectral measure $\nu$ has a strictly positive density, such that the function $\omega_{\nu}$ given by \eqref{omega:nu} is strictly positive and $(d + 1 + [\beta])$-times continuously differentiable on $\mathbb{S}^{d-1}$.
	\item[(H3)] The source term $f:[0,\infty)\to [0,\infty)$ fulfills the following conditions:
	\begin{itemize}
		\item[(a)] It is locally Lipschitz continuous, non-decreasing function with $f(0)=0$ and $f>0$ on $(0,\infty)$.
		\item[(b)] It satisfies the Osgood-type condition as defined in Definition \ref{Osgood-type:def} below.
	\end{itemize}
\end{itemize}
In particular, functions of Osgood-type satisfy the integral condition:
\[
\int_{1}^{\infty} \frac{ds}{f(s)} = \infty.
\]

The solution of \eqref{FE:0} can be represented using the {\it variation formula for Volterra equations} (cf. \cite{Pr93}) as follows
\begin{equation}
\label{IntegralSolution}
u(t,x)= S(t)u_0(x) + \int_{0}^{t} R(t-s)f(u(s,\cdot))(x)ds,
\end{equation}
where, for each $t>0$, $S(t)$ and $R(t)$ are linear and bounded operators on $L_p(\RR^d)$ of convolution type, i.e.,
\[
(S(t)v)(x) := \int_{\RR^d} Z(t,x-y)v(y)dy, \, \text{ and } \,  (R(t)v)(x):= \int_{\RR^d} Y(t,x-y) v(y)dy.
\] 
It is worth mentioning that the kernels $Z$ and $Y$ are generally not explicitly known. Let us provide a brief overview of the current understanding of the fundamental solution, represented by $Z(t,x)$, for \eqref{FE:0}. This solution is obtained when we set $u_0$ as the Dirac measure and $f$ as 0 in \eqref{FE:0}. For $\Psi_{\beta}(-i\nabla)=-\Delta$ with $\beta=2$ and $\omega_{\nu}\equiv 1$, and $\alpha\in (0,1)$, it is known (see, for example, \cite{ScWy89, Koch90}) that
\[
Z(t,x)=\pi^{-\frac{d}{2}}t^{\alpha-1}|x|^{-d} H^{20}_{12}\left(\frac{1}{4}|x|^2t^{-\alpha}\big|^{(\alpha,\alpha)}_{(d/2,1), (1,1)}\right),\quad t>0,\,x\in \RR^d\setminus\{0\},
\]
where $H$ represents the Fox $H$-function (\cite{KiSa04}, \cite{KiST06}). This representation of $Z$ is not very useful for deriving direct estimates for $Z$ due to the complexity of the $H$-function. However, by utilizing the analytic and asymptotic properties of $H$, the authors of \cite{EiKo04} (and also \cite{Koch90}) were able to derive sharp estimates for $Z$. Alternatively, in \cite{KeSVZ16}, they obtained $Z(t,x)$ by using the subordination principle for abstract Volterra equations with completely positive kernels, as detailed in Pr\"uss \cite[Chapter 4]{Pr93} and Cl\'ement-Nohel \cite{ClNo79} (and also \cite{PoVe18} and references therein). Specifically, $Z(t,x)$ is derived from the heat kernel $p_t(x) := \frac{1}{(4\pi t)^{d/2}}e^{-|x|^2/4t}$ as follows:
\begin{equation*}
	Z(t,x) = -\int_0^{\infty} p_{s}(x)\omega(t,ds),\quad t>0, x\in \RR^d,
\end{equation*}
where $-\omega(t,ds)$ is a probability measure on $\RR_+$ for each $t>0$.

The authors in \cite{KeSZ17} (also in \cite{KiLi16}) consider the case where $\Psi_{\beta}( -i\nabla)=(-\Delta)^{\beta/2}$ with $\beta\in (0,2)$ and $\alpha\in (0,1)$. They use the method from \cite{EiKo04} to construct the corresponding kernel $Z$ and estimate it. However, \cite[Section 2]{JoKo19} and \cite[Section 8.2]{Kolo09} provide a crucial development on this topic. The authors in these papers show that the linear Cauchy problem \eqref{FE:0} admits a fundamental solution $Z$, which is given by
\begin{equation}\label{Z:G}
	Z(t,x)=\dfrac{1}{\alpha}\displaystyle\int_0^\infty G(t^\alpha s,x)s^{-1-\frac{1}{\alpha}}G_\alpha(1,s^{-\frac{1}{\alpha}})ds.
\end{equation}
Here, $G$ is the Green function that solves the equation $\partial_t v(t,x)+\Psi_{\beta}(-i\nabla)v(t,x)=0$ with the initial condition $G(t,x)|_{t=0}=\delta_0(x)$, where $\delta_0$ is the Dirac delta distribution. Additionally, $G_{\alpha}(\cdot,\cdot)$ is the Green function that solves the problem $\partial_t v(t,s)+\dfrac{d^\alpha}{ds^\alpha}v(t,s)=0$ with the initial condition $G_{\alpha}(0,s)=\delta_0(s)$, where $\alpha\in (0,1)$ and $\dfrac{d^\alpha}{ds^\alpha}f(s):=\dfrac{1}{\Gamma(-\alpha)}\int_0^\infty \dfrac{f(s-\tau)-f(s)}{\tau^{1+\alpha}}d\tau$, as given in \cite[Formulas (1.111) and (2.74)]{Kolo19}. By combining \eqref{Z:G} with the approach of Aronson \cite{Aron67}, the authors in \cite{JoKo19} derive double-sided estimates for $Z(t,x)$ under the condition (H2) on $\nu$. It is worth noting that the representation of $Z$ given in \eqref{Z:G} has several advantages over the previously mentioned representations. Firstly, it is applicable to a wider range of operators $\Psi_{\beta}(-i\nabla)$, including fractional Laplacians of order $\beta\in(0,2)$. Secondly, it is expressed in terms of the Green function $G$, which has a more intuitive interpretation and is more widely studied in the literature. Thirdly, the representation is more amenable to analysis, as it can be used to obtain explicit estimates for $Z$ in certain cases, as demonstrated in \cite{JoKo19}. More precisely, for a fixed $T > 0$, $(t, x, y) \in (0, T] \times \RR^d \times \RR^d$ the following two-sided estimates for $Z(t,x-y)$ hold, for $\Omega:=|x-y|^{\beta} t^{-\alpha}\leq 1$ we have
\begin{equation}\label{z:doble:est}
Z(t,x-y) \asymp 
\begin{cases}
	t^{-\frac{d\alpha}{\beta}} & d<\beta,\\
	t^{-\alpha} (|\log \Omega| + 1) & d = \beta,\\
	t^{-\frac{d\alpha}{\beta}} \Omega^{1-\frac{d}{\beta}} & d > \beta,
\end{cases}
\end{equation}
and for $\Omega\geq 1$ the estimate
\[
Z(t,x-y) \asymp t^{-\frac{d\alpha}{\beta}} \Omega^{-1-\frac{d}{\beta}}.
\]
The notation $f(x)\asymp g(x)$ in $D$ above means that there exists constants $C, c > 0$ such that $f$ satisfies the following two-sided estimate, $cg(x) \leq f(x) \leq Cg(x)$, for all $x \in D$, for some region $D$.

Regarding the kernel $Y$, it is defined as the unique solution of the equation
\begin{equation}\label{Z:Y}
	Z(t,x) = (g_{1-\alpha}\ast Y(\cdot,x))(t),
\end{equation}
which also enjoys of two-sided estimates similar to \eqref{z:doble:est} as follows: for $\Omega\leq 1$,
\begin{equation}\label{y:doble:est}
Y(t,x-y) \asymp 
\begin{cases}
	t^{-\frac{d\alpha}{\beta}+\alpha-1} & d<2\beta,\\
	t^{-\alpha-1} (|\log \Omega| + 1) & d = 2\beta,\\
	t^{-\frac{d\alpha}{\beta}+\alpha-1} \Omega^{2-\frac{d}{\beta}} & d > 2\beta,
\end{cases}
\end{equation}
and for $\Omega\geq 1$ 
\[
Y(t,x-y) \asymp t^{-\frac{d\alpha}{\beta}+\alpha-1} \Omega^{-1-\frac{d}{\beta}},
\]
see \cite[Proposition 2.2]{SoVe22} and \cite[Lemma 2.15]{SoVe22} for the derivation of \eqref{Z:Y}.

\medbreak

 In this context, our work establishes conditions for the non-existence of local solutions to \eqref{FE:0} in the sense of \cite[p. 78]{QS07}. This demonstrates that the Osgood-type condition of the non-linear term does not ensure the existence of positive solutions. Specifically, we demonstrate the occurrence of \textit{instantaneous blow-up}, a phenomenon not previously observed in the non-Gaussian case to the best of our knowledge. While this phenomenon has been studied extensively in other setting, such as in the case of $\beta=2$, $\alpha=1$ with $\omega_{\nu}(\theta) \equiv 1$ investigated in \cite{LaRoSi13}, and for $\beta=2$, $\alpha \in (0,1)$ with $\omega_{\nu}(\theta) \equiv 1$ explored in \cite{ZS15}, as well as global existence and \textit{blow-up in finite time} scenarios, with $\beta\in (0,2)$, $\alpha \in (0,1)$ and $\omega_{\nu}(\theta)$ as in (H2), discussed in \cite{SoVe22} and \cite{SoVe23}, respectively. Nevertheless, we also demonstrate a scenario where a global solution is attainable even with an Osgood-type non-linear term.

A key strength of the present study is the use of \textit{critical exponents} for non-existence and existence of solutions to \eqref{FE:0}, respectively. From the literature we know that these exponents are used to express inequalities that lead to well-posed or ill-posed problems. In both cases, the definition of critical exponents is required in this work. In Section \ref{sec:3}, for instance, $q_c:=k\left(1-\frac{\beta}{\alpha d+\beta}\right)$ is considered the critical exponent leading to an ill-posed problem whenever $q<q_c$, while in the last section another critical exponent is defined for ensuring the existence of global solutions.

\medbreak
 
This paper is organized as follows. Section \ref{sec:2} is devoted to review some properties of the fundamental solutions $Z,Y$ of problem \eqref{FE:0}, as well as we prove new properties of them. The main results of this section are Lemma \ref{LemmaL_1} and Theorem \ref{GlobalSolution_boundedu_0}. In Section \ref{sec:3} we prove our main result for instantaneous blow-up, in the non-Gaussian context, which is stated in Theorem \ref{blow-up}. The final section shows a result for global solutions and the critical exponent is introduced in terms of the parameters of the non-Gaussian equation \eqref{FE:0}. The main result of this section is given by Theorem \ref{globalsolutionsmallu0}.
\section{Preliminaries}
\label{sec:2}
In this section, we show some additional properties of the kernels $Z,Y$ for proving our instantaneous blow-up result. Furthermore, we revisit a key construction technique for Osgood-type functions, as introduced in \cite[Section 3]{LaRoSi13}, laying the groundwork for our subsequent discussions.

In addition, the study of Gaussian evolution equations has gave inspiration for using another approach to existence of solutions. For instance, in \cite{RoSi13} the authors employ \textit{super-solution} methods.  In order to do an extension to the non-Gaussian case, let us define the non-linear operator
\begin{equation}\label{F:u:S:R:f}
\mathcal{F}u(t):= S(t)u_0 + \int_{0}^{t} R(t-s)f(u(s,\cdot))ds,
\end{equation}
for measurable functions $u :[0, T)\times \RR^d \rightarrow [0, \infty]$.

It is worth mentioning that in the reaction-diffusion equation, the operators $S$ and $R$ coincide, i.e., 
\[ 
S(t)v = R(t)v = \frac{1}{(4\pi t)^{d/2}}\int_{\mathbb{R}^d} e^{-|x-y|^2/4t}v(y) dy, 
\]
since in this case the kernels $Z$ and $Y$ match with the heat kernel. For the reader's convenience, we compile certain properties of $Z$ and $Y$ from \cite{SoVe22}, which we later utilize throughout the paper.

\begin{proposition}\label{propo:Z:Y}
	Let $\alpha\in (0,1)$ and $\beta\in (0,2)$. The following properties hold:
	\begin{itemize}
		\item[(i)] $\int_{\RR^d} Z(t,x)dx =1$ for all $t>0$, and $Z(t,x) = t^{-\alpha d/ \beta} Z(1, t^{-\alpha/\beta} x)$ for all $t>0$ and $x\in \RR^d$.
		\item[(ii)] $\int_{\RR^d} Y(t,x)dx =g_{\alpha}(t)$ for all $t>0$, and $Y(t,x) = t^{-\alpha d/ \beta + \alpha -1} Y(1, t^{-\alpha/\beta} x)$ for all $t>0$ and $x\in \RR^d$.
	\end{itemize}
\end{proposition}
For the proof, refer to \cite[Lemma 2.12]{SoVe22}. The first part of (ii) can be found in the proof of Theorem 2.14 in \cite{SoVe22}. The $L_p$-properties of $Z$ and $Y$ are gathered from \cite[Theorem 2.8 and Theorem 2.10]{SoVe22}, respectively, in the following proposition.

\begin{proposition}\label{propo:lp:Z:Y}
		Let $d\in \NN$, $\alpha\in (0,1)$, and $\beta\in (0,2)$. Under assumption $(H2)$, the following properties are established:
		\begin{itemize}
			\item[(i)] The kernel $Z(t,\cdot)$ is in $L_p(\mathbb{R}^d)$ for all $t>0$ if and only if $1\leq p<\kappa_1$, where
			\[
			\kappa_1=\kappa_1(d,\beta):=
			\begin{cases}
				\frac{d}{d-\beta} & \text{if } d>\beta,\\
				\infty & \text{otherwise}.
			\end{cases}
			\]
			Moreover, the two-sided estimate 
			\begin{equation}
				\label{cotasZp}
				\norm{Z(t,\cdot)}_p\asymp t^{-\frac{\alpha d}{\beta}\left(1-\frac{1}{p}\right)},\; t>0
			\end{equation}
			holds for every $1\leq p<\kappa_1$. In the case of $d<\beta$, \eqref{cotasZp} remains true for $p=\infty$.\\
			\item[(ii)] The kernel $Y(t,\cdot)$ is in $L_p(\mathbb{R}^d)$ for all $t>0$ if and only if $1\leq p<\kappa_2$, where
			\[
			\kappa_2=\kappa_2(d,\beta):=
			\begin{cases}
				\frac{d}{d-2\beta} & \text{if } d>2\beta,\\
				\infty & \text{otherwise}.
			\end{cases}
			\]
			Moreover, the two-sided estimate 
			\begin{equation}
				\label{cotasYp}
				\norm{Y(t,\cdot)}_p\asymp t^{-\frac{\alpha d}{\beta}\left(1-\frac{1}{p}\right)+(\alpha-1)},\; t>0
			\end{equation}
			holds for every $1\leq p<\kappa_2$. In the case of $d<2\beta$, \eqref{cotasYp} remains true for $p=\infty$.
		\end{itemize}
\end{proposition}

\begin{definition}
\label{DefSolution}
Let $u_0$, $f$ be non-negative functions, and let $\mathcal{F}$ be given by \eqref{F:u:S:R:f}. The function $u :[0, T)\times \RR^d  \rightarrow [0, \infty]$ is an {\bf integral solution} of \eqref{FE:0} if it is measurable, finite almost everywhere, and satisfies $u=\mathcal{F}u$ for almost every $(t, x)\in  (0, T)\times \RR^d$. If $T$ can be arbitrarily large, we say that the solution is {\bf global}.
\end{definition}
That is, a solution in the sense of Definition \ref{DefSolution} satisfies \eqref{IntegralSolution}. 
\begin{definition}
\label{DefSuperSolution}
Let $u_0$, $f$ be non-negative functions, and let $\mathcal{F}$ be given by \eqref{F:u:S:R:f}. A function $u :[0, T)\times \RR^d  \rightarrow [0, \infty]$ is a {\bf super-solution} ({\bf sub-solution}) of \eqref{FE:0} if it is measurable, finite almost everywhere, and satisfies $u\geq\mathcal{F}u$ ($u\leq\mathcal{F}u$) for almost every $(t, x)\in  (0, T)\times \RR^d$.
\end{definition}
\begin{lemma}
\label{LemaSuperSol}
Let $f:[0,\infty)\rightarrow [0,\infty)$ be a continuous and non-decreasing function. Suppose $u_0\in L_q(\RR^d)$ is non-negative, where $1\leq q<\infty$. Then the operator $\mathcal{F}$ defined by \eqref{F:u:S:R:f} admits a solution in $(0, T)\times \RR^d$ if and only if it admits a super-solution in $(0, T)\times \RR^d$.
\end{lemma}
\begin{proof}
It is clear that every solution is a super-solution. Now, let $u$ be a super-solution of $\mathcal{F}$, that is $u\geq\mathcal{F}u$. Firstly, we note that $\mathcal{F}$ is a monotonic operator. Indeed, suppose $u_1\geq u_2.$ It follows that
\begin{align*}
\left(\mathcal{F}u_1-\mathcal{F}u_2\right)(t,x)&=\int_{0}^{t} R(t-s)\left(f(u_1(s,\cdot))-f(u_2(s,\cdot))\right)(x)ds\\
&=\int_0^t \int_{\RR^d}Y(t-s,x-y)\left(f(u_1(s,y))-f(u_2(s,y))\right)dyds\geq 0
\end{align*}
a.e. in $(0, T)\times \RR^d$, due to the positivity of kernel $Y$ and the monotonicity of $f$.

Consequently, we have that $\mathcal{F}u\geq\mathcal{F}^2u$ and inductively we can define a decreasing sequence of functions 
\[
\mathcal{F}u\geq\mathcal{F}^2u\geq\cdots \geq\mathcal{F}^{n}u\geq\mathcal{F}^{n+1}u,\quad n\geq 1,
\]
where $0\leq u_n:=\mathcal{F}^n(u)\leq u$ for all $n$. Therefore, we set
\[
v(t,x)=\lim_{n\rightarrow\infty} u_n(t,x)
\]
a.e. in $(0, T)\times \RR^d$. Next, we show that $v$ satisfies \eqref{IntegralSolution}.
\begin{align*}
u_{n+1}&=\mathcal{F}u_n(t,x)\\
&= S(t)u_0(x) + \int_{0}^{t} R(t-s)f(u_n(s,\cdot))(x)ds\\
&=\int_{\RR^d}Z(t,x-y)u_0(y)dy + \int_0^t \int_{\RR^d}Y(t-s,x-y)f(u_n(s,y))dyds.
\end{align*}
Utilizing the monotone convergence theorem together with the continuity of $f$, yields
\[
\lim_{n\rightarrow\infty}u_{n+1}=\int_{\RR^d}Z(t,x-y)u_0(y)dy + \int_0^t \int_{\RR^d}Y(t-s,x-y)f\left(\lim_{n\rightarrow\infty}u_n(s,y)\right)dyds.
\]
This shows that $v=\mathcal{F}v$.
\end{proof}
\begin{proposition}{\label{propo:z:phi}}
	Fix $\tau\in (0, d)$ and $R>1$ and let $u_0\in L_1(\RR^d)$ be the non-negative, radially symmetric function given by
	\[
	u_0(x) := |x|^{-\tau}1_{B(0,R)}(x) :=
	\begin{cases}
		|x|^{-\tau}, & |x|\leq R,\\
		0, & |x| >R.
	\end{cases}
	\]
	Let $z(t,x):=(Z(t, \cdot)\star u_0)(x)$ and $M:=\min\left\{z(t,\hat{x}) \, : \,  \hat{x}\in \mathbb{S}^{d-1}, \, 0\leq t\leq 1\right\}$. If $\rho\in (0,\alpha/\beta)$ then for any $\phi\geq M$ we have
	\[
	z(t,x) \geq \phi \, \text{ for } t^{\alpha/\beta} \leq |x| \leq t^{\rho}, \text{ for all } 0<t\leq \left(\frac{\phi}{M}\right)^{-\frac{1}{\tau \rho}}.
	\]
\end{proposition}
\begin{proof}
For any $0<t\leq 1$ and $x=|x|\hat{x}$, $\hat{x}\in \mathbb{S}^{d-1}$, we have
	\begin{align*}
		z(t, x)& = \int_{\RR^d} Z(t, x - y) u_0(y)dy\\
		& =|x|^d\int_{\RR^d} Z(t, |x|\hat{x} -  |x|w) u_0(|x| w)dw\\
		& =|x|^d\left( \int_{\RR^d \setminus B(0, |x|^{-1}R)} Z(t, |x|\hat{x} -  |x|w) u_0(|x| w)dw +\int_{B(0, |x|^{-1}R)} Z(t, |x|\hat{x} -  |x|w) u_0(|x| w)dw\right)\\
		& = |x|^{d-\tau}\int_{B(0, |x|^{-1}R)} Z(t, |x|(\hat{x} -  w))|w|^{-\tau} dw.
\end{align*}		
		Now, by the scaling property of $Z$ in Proposition \ref{propo:Z:Y}(i), we obtain
\begin{align*}
z(t, x) & =|x|^{d-\tau}\int_{B(0, |x|^{-1}R)} t^{-\alpha d/ \beta} Z(1, t^{-\alpha / \beta} |x|(\hat{x} -  w))|w|^{-\tau} dw\\
& =|x|^{d-\tau}\int_{B(0, |x|^{-1}R)} t^{-\alpha d/ \beta} Z(1,(t |x|^{-\beta/\alpha})^{-\alpha / \beta}(\hat{x} -  w))|w|^{-\tau} dw\\
& =|x|^{-\tau}\int_{B(0, |x|^{-1}R)} Z(t |x|^{-\beta/\alpha},\hat{x} -  w)|w|^{-\tau} dw.
\end{align*}
By hypothesis it follows that $|x|\leq 1$ and that $t |x|^{-\beta/\alpha}\leq 1$, therefore
\begin{align*}
z(t, x) & \geq |x|^{-\tau}\int_{B(0, R)} Z(t |x|^{-\beta/\alpha},\hat{x} -  w)|w|^{-\tau} dw\\
&= |x|^{-\tau}\int_{\RR^d} Z(t |x|^{-\beta/\alpha},\hat{x} -  w)u_0(w) dw\\
&=|x|^{-\tau}z(t |x|^{-\beta/\alpha},\hat{x})\\
&\geq |x|^{-\tau}M\geq t^{-\tau\rho}M\geq \phi.
\end{align*}
\end{proof}

Next, we present a family of Osgood-type functions $(f^{[k]})_{k>1}$ whose construction is detailed in \cite[Section 3]{LaRoSi13}. However, for the convenience of the reader, we summarize its main characteristics.

Consider $k>1$, and fix a number $\phi_0>2^{\frac{1}{k-1}}$. Define the sequence $(\phi_i)_{i\in \NN}$ as follows
\[
\phi_{i}=\phi_{i-1}^k,\quad i\geq 1,
\]
ensuring
\[
1<\phi_{i-1}<\frac{\phi_i}{2},\quad i\geq 1, 
\]
and 
$\phi_i\rightarrow\infty$ as $i\rightarrow\infty$.

Now, we define the function $f^{[k]}$ as
\begin{equation}
\label{Osgoodf}
f^{[k]}(s) :=\begin{cases}\left(1-\phi_0^{1-k}\right)s^k & \text{if }s\in J_0=[0,\phi_0],\\ \phi_i-\phi_{i-1} & \text{if }s\in I_i=\left[\phi_{i-1},\frac{\phi_i}{2}\right],~i\geq 1,\\l_i(s) &\text{if }s\in J_i=\left(\frac{\phi_i}{2},\phi_i\right),~i\geq 1,  \end{cases}
\end{equation}
where $l_i$ is a polynomial of order $1$ that interpolates between the values of $f^{[k]}$ at $\frac{\phi_i}{2}$ and $\phi_i$. It is important to note that each $f^{[k]}$ satisfies assumptions (H3)(a) and the integral condition:
\[
\int_{1}^{\infty} \frac{ds}{f^{[k]}(s)} = \infty,
\]
as demonstrated in \cite[Section 3]{LaRoSi13}.

It follows that the function $\tilde{f}^{[k]} : [0,\infty)\rightarrow [0,\infty)$, given by
\begin{equation}
\label{ftilde}
\tilde{f}^{[k]}(s) :=\begin{cases}0 & \text{if }s\in J_0=[0,\phi_0],\\ \phi_i-\phi_{i-1} & \text{if }s\in I_i\cup J_1=\left[\phi_{i-1},\phi_i\right),~i\geq 1,\end{cases}
\end{equation}
satisfies $f^{[k]}\geq  \tilde{f}^{[k]}$ on $[0,\infty)$.

\begin{definition}\label{Osgood-type:def}
	We define a continuous function $f:[0,\infty)\to [0,\infty)$ to be of {\bf Osgood-type} if there exists $k>1$ such that $f^{[k]}$, as defined in \eqref{Osgoodf}, satisfies the two-sided estimate $f(s)\asymp f^{[k]}(s)$ for all $s\geq 0$.
\end{definition}

\begin{lemma}
	\label{LemmaL_1}
	Suppose that $f$ satisfies condition (H3) for some $k>1+\frac{\beta}{\alpha d}$. Then there exists a non-negative function $u_0\in L_1(\RR^d)$ such that there is no local integral solution $u(\cdot)$ of \eqref{FE:0} that remains in $L_1(\RR^d)$ for all $t>0$.
\end{lemma}
\begin{proof}
	We set $0<\tau<d$ and define $u_0$ as in Proposition \ref{propo:z:phi}. Suppose, by contradiction, that there exists a measurable function $u:[0,T)\times \mathbb{R}^d\rightarrow [0,\infty]$, finite almost everywhere, which is a local integral solution of \eqref{FE:0} with $u(0)=u_0\geq 0$. Since $f$ is non-negative and non-decreasing, it follows that
	\begin{align*}
		u(t,x) & = S(t)u_0(x) + \int_{0}^{t} R(t-s)f(u(s,\cdot))(x)ds\\
		& \geq S(t) u_0(x),
	\end{align*}
and thus,
\[
f(u(t,x)) \geq f(S(t)u_0(x)).
\]
Therefore, we estimate the $L_1$-norm of $u$ by utilizing Fubini's Theorem and the integral property of $Y$ as stated in Proposition \ref{propo:Z:Y}(ii). We obtain:
\begin{align*}
	|u(t,\cdot)|_1 & \geq \int_{\RR^d} \int_{0}^{t} R(t-s)f(S(s)u_0)(x)dsdx\\
	& = \int_{\RR^d} \int_{0}^{t} \int_{\RR^d} Y(t-s,x-y) f(S(s)u_0(y))dy ds dx\\
	& = \int_{0}^{t} g_{\alpha}(t-s) \left(\int_{\RR^d} f(S(s)u_0(y))dy\right)ds.
\end{align*}

Now, letting  $z(t,x)=(Z(t, \cdot)\star u_0)(x) = S(t)u_0(x)$, using the function $\tilde{f}:=\tilde{f}^{[k]}$ introduced in \eqref{ftilde}, together with the fact $f\geq c\tilde{f}$ and $\phi_i\leq\frac{\phi_{i+1}}{2}$, we estimate
\begin{align*}
	\int_{\RR^d} f(S(s)u_0(y))dy  & \gtrsim\int_{\RR^d}\tilde{f}(S(s)u_0(y))dy\\
	& =\int_{\RR^d}\tilde{f}(z(s,y))dy\\
	& \gtrsim\sum_{i=N}^{\infty} \left|\left\{y \, : \, \phi_i\leq z(s,y) < \phi_{i+1} \right\}\right|(\phi_{i+1} - \phi_i)\\
	& \gtrsim\sum_{i=N}^{\infty} \left|\left\{y \, : \,  \phi_i \leq z(s,y) <\phi_{i+1} \right\}\right| \frac{\phi_{i+1}}{2}\\
	& =\sum_{i=N}^{\infty} \left|\left\{y \, : \, \phi_i \leq z(s,y) < \phi_{i+1} \right\}\right|\frac{\phi_{i}^k}{2}\\
	&\gtrsim \frac{\phi_{N}^k}{2}\sum_{i=N}^{\infty} \left|\left\{y \, : \, \phi_i \leq z(s,y) < \phi_{i+1} \right\}\right|\\
	&= \frac{\phi_{N}^k}{2}\left|\left\{y \, : \,  z(s,y)>\phi_{N} \right\}\right|
\end{align*}
for some $N\in\mathbb{N}$.

Returning to the estimate of the $L_1$-norm of $u$ and using this last inequality, taking $i$ sufficiently large so that $\phi_i\geq M$ and $(\phi_i/M)^{-\frac{1}{\tau \rho}} \leq t$, it follows from Proposition \ref{propo:z:phi} that
\begin{align*}
	|u(t,\cdot)|_1 & \gtrsim\int_{0}^{t} g_{\alpha}(t-s) \left|\left\{x \, : \, z(s,x) > \phi_i \right\}\right|\frac{\phi_{i}^k}{2} ds\\
	& \gtrsim\int_{0}^{(\phi_i/M)^{-\frac{1}{\tau \rho}}} g_{\alpha}(t-s) \left|\left\{x \, : \, z(s,x) > \phi_i \right\}\right|\frac{\phi_{i}^k}{2} ds\\
	& \gtrsim w_d \frac{\phi_{i}^k}{2}  \int_{0}^{(\phi_i/M)^{-\frac{1}{\tau \rho}}} g_{\alpha}(t-s)(s^{d\rho}-s^{\alpha d/\beta})ds\\
	& = w_d \frac{\phi_{i}^k}{2}  \int_{0}^{(\phi_i/M)^{-\frac{1}{\tau \rho}}}\frac{(t-s)^{\alpha-1}}{\Gamma(\alpha)} (s^{d\rho}-s^{\alpha d/\beta})ds\\
	& \gtrsim w_d \frac{t^{\alpha-1}}{\Gamma(\alpha)} \frac{\phi_{i}^k}{2}  \int_{0}^{(\phi_i/M)^{-\frac{1}{\tau \rho}}}(s^{d\rho}-s^{\alpha d/\beta})ds\\
	& =C(d, \alpha, c)~t^{\alpha-1} \left(\phi_{i}^{k - \frac{1}{\tau \rho}(d\rho + 1)}-\phi_{i}^{k - \frac{1}{\tau \rho}(\alpha d/\beta + 1)}\right),
\end{align*}
where $w_d$ denotes the volume of the unit ball in $\RR^d$. In this way, if $k>1+\frac{\beta}{\alpha d}$, then it is possible to choose $\tau\in (0,d)$ and $\rho\in\left(0,\frac{\alpha}{\beta}\right)$ such that $k > \frac{1}{\tau}\left(d + \frac{1}{\rho}\right)$ and 
\[
\phi_{i}^{k - \frac{1}{\tau \rho}(d\rho + 1)}-\phi_{i}^{k - \frac{1}{\tau \rho}(\alpha d/\beta + 1)} \to \infty \text{ as } i\to \infty.
\]
Hence, for this choice of $u_0\in L_1(\RR^d)$ there is no integral solution of \eqref{FE:0} that remains in $L_1(\RR^d)$.
\end{proof} 

Finally, we show that if $u_0$ is bounded then the solution exists globally in time.

\begin{theorem}
\label{GlobalSolution_boundedu_0}
Suppose that $f$ satisfies condition (H3). If $u_0\in L_\infty (\RR^d)$ then the problem \eqref{FE:0} has a unique global solution $u\in L_{\infty,loc}((0,\infty);L_\infty(\RR^d))$ in the sense of Definition \ref{DefSolution}.
\end{theorem}
\begin{proof}
We define the operator
\[
\mathcal{M}v(t,x):=\int_{\RR^d}Z(t,x-y)v_0(y)dy + \int_0^t \int_{\RR^d}Y(t-s,x-y)g(v(s,y))dyds
\]
on the Banach space $L_{\infty,loc}((0,\infty);L_\infty(\RR^d))$, where $g$ is a non-decreasing Lipschitz function with $g(0) = 0$ and $v_0\in L_\infty(\RR^d)$. As in the proof of \cite[Lemma 1.3]{AE87}, using estimates of $Z,Y$ given in \eqref{z:doble:est} and \eqref{y:doble:est} respectively, we derive that the operator $\mathcal{M}$ has a unique fixed point $v$. Furthermore, $v\geq w$ whenever $v_0\geq w_0$, where $w$ is the fixed point associated with $w_0\in L_\infty(\RR^d)$. We apply this result to the sequence 
\[
f_n(s):=\begin{cases} 0 &\text{if } s<0,\\ f(s) &\text{if } s\in J_0\cup_{i=1}^n I_i\cup J_i,\\ f(\phi_n) &\text{if } s\geq\phi_n,\\ 
\end{cases}
\]
for $n\geq 1$. Observe that the local Lipschitz continuity of $f$ and the compactness of the interval $[0,\phi_n]$ ensure that each $f_n$ is Lipschitz continuous. In this way, there exists a unique function $u_n\in L_{\infty,loc}((0,\infty);L_\infty(\RR^d))$ such that $0\leq u_n$ and
\[
u_n(t,x)=\int_{\RR^d}Z(t,x-y)\left(u_0+\frac{1}{n}\right)(y)dy + \int_0^t \int_{\RR^d}Y(t-s,x-y)f_n(u_n(s,y))dyds,
\]
for $x\in\RR^d$ and $0<t<\infty$. Since $\frac{1}{n}\geq\frac{1}{n+1}$ we have that $u_{n+1}\leq u_{n}$ for all $n\geq 1$. Thus, for almost every $(t,x)\in (0,\infty)\times\RR^d$, the sequence of real numbers $(u_n(t,x))_{n\in\mathbb{N}}$ is decreasing and bounded from below by zero. Consequently, we can define the function
\[
u(t,x)=\lim_{n\rightarrow\infty}u_n(t,x)
\]
a.e. in $(0,\infty)\times\RR^d$. Since $u_1\in L_{\infty,loc}((0,\infty);L_\infty(\RR^d))$, for any $T>0$ there exists a constant $C(T)$ such that $u_n(t,x)\leq C(T)$ for all $n\geq 1$ and $(t,x)\in (0,T)\times\RR^d$. On the other hand, for $N\in\mathbb{N}$ sufficiently large we have $C(T)<\phi_ N$ and hence $f_n(u_n(t,x))=f(u_n(t,x))$ whenever $n\geq N$. Therefore,
\[
u_n(t,x)=\int_{\RR^d}Z(t,x-y)\left(u_0+\frac{1}{n}\right)(y)dy + \int_0^t \int_{\RR^d}Y(t-s,x-y)f(u_n(s,y))dyds
\]
for all $n\geq N$ and $(t,x)\in (0,T)\times\RR^d$. The continuity of $f$ and the monotonic convergence theorem yield
\[
u(t,x)=\int_{\RR^d}Z(t,x-y)u_0 (y)dy + \int_0^t \int_{\RR^d}Y(t-s,x-y)f(u(s,y))dyds
\]
as $n\rightarrow\infty$.

For the uniqueness we suppose that $\widetilde{u}$ is also an integral solution to problem \eqref{FE:0}. Since each $f_n$ is Lipschitz continuous with Lipschitz constant $C_n$, we can find $N\in\mathbb{N}$ as above, such that 
\begin{align*}
|u(t,x)-\widetilde{u}(t,x)|&\leq\int_0^t \int_{\RR^d}Y(t-s,x-y)|f(u(s,y))-f(\widetilde{u}(s,y))|dyds\\
&=\int_0^t \int_{\RR^d}Y(t-s,x-y)|f_N(u(s,y))-f_N(\widetilde{u}(s,y))|dyds\\
&\leq C_N\int_0^t \int_{\RR^d}Y(t-s,x-y)|u(s,y)-\widetilde{u}(s,y)|dyds\\
&\leq C_N\int_0^t \int_{\RR^d}Y(t-s,x-y)\norm{u(s)-\widetilde{u}(s)}_\infty dyds\\
&\leq C(N,\alpha)\int_0^t (t-s)^{\alpha-1}\norm{u(s)-\widetilde{u}(s)}_\infty ds
\end{align*}
for all $(t,x)\in (0,T)\times\RR^d$. Hence, Gronwall’s inequality (see \cite[Corollary 2]{YGD07}) shows that $u(t)=\widetilde{u}(t)$.
\end{proof}
 
\section{Instantaneous blow-up}
\label{sec:3}
In this section, we demonstrate that problem \eqref{FE:0} does not admit a local solution in the sense of Definition \ref{DefSolution}. To achieve this, we establish that the function $u(t)$ is not locally integrable in $\mathbb{R}^d$ for all $t>0$, and consequently it does not belong to $L_q(\mathbb{R}^d)$ for any $1\leq q <\infty$.

\begin{theorem}
\label{blow-up}
	Let $q\in [1,\infty)$, and suppose that $f$ satisfies condition (H3) for some $k > q\left(1+\frac{\beta}{\alpha d}\right)$.  Then there exists a function $u_0\in L_q(\mathbb{R}^d)$ such that \eqref{FE:0} possesses no local integral solution. Moreover, any solution $u(\cdot)$ that satisfies \eqref{IntegralSolution} is not in $L_{1,\text{loc}}(\mathbb{R}^d)$ for any $t > 0$.
\end{theorem}
\begin{proof}
Similar to Proposition \ref{propo:z:phi}, we define $u_0$ as 
\[
	u_0(x) := |x|^{-\tau}1_{B(0,R)}(x) :=
	\begin{cases}
		|x|^{-\tau}, & |x|\leq R,\\
		0, & |x| >R,
	\end{cases}
	\]
with $\tau q<d$. Let $t\in (0,1]$ and take $i$ sufficiently large so that $(\phi_i/M)^{-\frac{1}{\tau \rho}} \leq t$. We denote by $B(\epsilon)$ the closed ball in $\RR^d$, with center at the origin and radius $\epsilon >1$.

The construction of $f^{[k]}$ shows that $f\gtrsim \frac{\phi_{i+1}}{2}$ on $[\phi_i,\infty)$, $i\geq 1$, and similar arguments in the proof of Lemma \ref{LemmaL_1} yield
\begin{align*}
	\int_{B(\epsilon)}u(t,x)dx & \geq \int_{B(\epsilon)} \int_{0}^{t} R(t-s)f(S(s)u_0)(x)ds\,dx\\
	& = \int_{0}^{t} \int_{B(\epsilon)} \int_{\RR^d} Y(t-s,x-y) f(z(s,y)) dy\,dx\,ds\\
    & \geq \int_{0}^{(\phi_i/M)^{-\frac{1}{\tau \rho}}}  \int_{s^{\alpha/\beta}\leq |y|\leq s^{\rho}} \int_{B(\epsilon)} Y(t-s,x-y) f(z(s,y))dx\,dy\,ds\\	
	& \gtrsim \int_{0}^{(\phi_i/M)^{-\frac{1}{\tau \rho}}}  \int_{s^{\alpha/\beta}\leq |y|\leq s^{\rho}} \int_{B(\epsilon)} Y(t-s,x-y)\frac{\phi_{i+1}}{2} dx\,dy\,ds\\
	& = \frac{\phi_{i+1}}{2}\int_{0}^{(\phi_i/M)^{-\frac{1}{\tau \rho}}}  \int_{s^{\alpha/\beta}\leq |y|\leq s^{\rho}} \int_{|x+y|\leq\epsilon} Y(t-s,x)dx\,dy\,ds.
\end{align*}
Now, since $|y|\leq 1$ we note that if $|x|\leq\epsilon-1$ then $|x+y|\leq \epsilon$. Therefore
\[
	\int_{|x+y|\leq\epsilon} Y(t-s,x)dx\geq \int_{|x|\leq\epsilon-1} Y(t-s,x)dx.
\]
Let $\delta=\min(\epsilon-1,1)$. The scaling property of $Y$ as stated in Proposition \ref{propo:Z:Y}(ii), along with estimates \eqref{y:doble:est} and the fact that $t\leq 1$, yield:
\begin{align*}
\int_{|x|\leq\epsilon-1} Y(t-s,x)dx &=(t-s)^{\alpha-1}\int_{|x|\leq \epsilon-1} Y(1,(t-s)^{-\alpha/\beta}x)(t-s)^{-\alpha d/\beta}dx\\
& =(t-s)^{\alpha-1}\int_{|x|\leq (t-s)^{-\alpha/\beta}(\epsilon-1)} Y(1,x)dx\\
&\geq t^{\alpha-1}\int_{|x|\leq \epsilon-1} Y(1,x)dx\\
&\geq t^{\alpha-1}\int_{|x|\leq \delta} Y(1,x)dx\\
&\geq C(d,\epsilon)t^{\alpha-1}.
\end{align*}
This implies that
\begin{align*}
	\int_{B(\epsilon)}u(t,x)dx & \geq C(d,\epsilon,c)t^{\alpha-1} \frac{\phi_{i+1}}{2}\int_{0}^{(\phi_i/M)^{-\frac{1}{\tau \rho}}}  \int_{s^{\alpha/\beta}\leq |y|\leq s^{\rho}} dy\,ds\\
	&= \frac{C(d,\epsilon,c)}{2} t^{\alpha-1}\phi_i^k\int_{0}^{(\phi_i/M)^{-\frac{1}{\tau \rho}}}(s^{\rho d}-s^{\alpha d/\beta})ds\\
	&= C(d,\epsilon,\alpha,c) t^{\alpha-1}\left(\phi_{i}^{k - \frac{1}{\tau \rho}(d\rho + 1)}-\phi_{i}^{k - \frac{1}{\tau \rho}(\alpha d/\beta + 1)}\right).
\end{align*}
Thus, if $k>q\left(1+\frac{\beta}{\alpha d}\right)$ then it is possible to choose $\tau \in\left(0,\frac{d}{q}\right)$ and $\rho\in\left(0,\frac{\alpha}{\beta}\right)$ such that $k > \frac{1}{\tau}\left(d + \frac{1}{\rho}\right)$ and 
\[
\int_{B(\epsilon)}u(t,x)dx \rightarrow \infty
\]
whenever $i\rightarrow \infty$. This shows that the blow-up is local to the origin. For $t>1$, we only note that the integral over $[0,t]$ is bigger than the integral over $[0,1]$. The proof is complete.
\end{proof}

\section{Global solutions}
\label{sec:4}
We establish global solutions in $L_q(\mathbb{R}^d)$ when the source term $f$ is an Osgood-type function and $u_0\in L_{q'}(\mathbb{R}^d)$, with $1\leq q<\infty$ and $q'$ being the \textit{critical exponent}.

\begin{theorem}
	\label{globalsolutionsmallu0}
	Assuming hypotheses \text{(H1)-(H3)} hold. Define $q':=\frac{d}{\beta}(k-1)\geq 1$, and let $u_0\in L_{q'}(\mathbb{R}^d)$ be a non-negative function. Suppose that $\max\left(\frac{d}{\beta},k,\frac{d(k-1)}{\beta}\right)<q<\frac{\alpha dk(k-1)}{\beta}$. If $|u_0|_{q'}$ is sufficiently small, then the problem \eqref{FE:0} has a global solution $u\in L_{\infty,\text{loc}}((0,\infty);L_q(\mathbb{R}^d))$ in the sense of Definition \ref{DefSolution}.
\end{theorem}
\begin{proof}
We define the Banach space
\[
E:=L_{\infty,loc}((0,\infty);L_q(\RR^d))
\]
with the norm
	\[
	\left\|v\right\|_{E}:=\esssup_{t> 0}t^{\frac{\alpha d}{\beta}\left(\frac{1}{q'}-\frac{1}{q}\right)}\norm{v(t,\cdot)}_{q}.
	\]
Let $\lambda >0$. On this space, we define the non-linear operator
\[
\mathcal{G}u(t)= S(t)u_0 + \lambda\int_{0}^{t} R(t-s)|u|^{k-1}u(s)ds.
\] 
Using the $L_p$-properties of kernels $Z$ and $Y$ as are stated in Proposition \ref{propo:lp:Z:Y}, together with Young's convolution inequality applied to the relationships $1+ \frac{1}{q} = \frac{1}{q'} + \frac{1}{p}$ and $1+ \frac{1}{q} = \frac{k}{q} + \frac{1}{p}$ respectively, we show that this operator is well-defined. Indeed, for $u\in E$ we have that
\begin{align*}
\norm{\mathcal{G}u(t)}_q & \leq \norm{S(t)u_0}_q + \lambda\int_{0}^{t} \norm{R(t-s)|u|^{k-1}u(s)}_q ds\\
&= \norm{Z(t)\star u_0}_q + \lambda\int_{0}^{t} \norm{Y(t-s)\star|u|^{k-1}u(s)}_q ds\\
&\lesssim t^{-\frac{\alpha d}{\beta}\left(\frac{1}{q'}-\frac{1}{q}\right)}\norm{u_0}_{q'}+
\lambda\int_{0}^{t} (t-s)^{-\frac{\alpha d}{\beta}\left(\frac{k}{q}-\frac{1}{q}\right)+\alpha-1}\norm{|u|^{k-1}u(s)}_{\frac{q}{k}} ds.
\end{align*}
We note that $\frac{k-1}{k}+\frac{1}{k}=1$. Then, applying Hölder's inequality in the second integral, we obtain
\[
\norm{|u|^{k-1}u(s)}_{\frac{q}{k}}\leq \norm{u(s)}_q^{k-1}\norm{u(s)}_q=\norm{u(s)}_q^k.
\]
Therefore, 
\begin{align*}
\norm{\mathcal{G}u(t)}_q&\lesssim t^{-\frac{\alpha d}{\beta}\left(\frac{1}{q'}-\frac{1}{q}\right)}\norm{u_0}_{q'}+
\lambda\left\|u\right\|_E^k\int_{0}^{t} (t-s)^{-\frac{\alpha d(k-1)}{\beta q}+\alpha-1}s^{-\frac{\alpha d k}{\beta}\left(\frac{1}{q'}-\frac{1}{q}\right)}ds\\
&=t^{-\frac{\alpha d}{\beta}\left(\frac{1}{q'}-\frac{1}{q}\right)}\norm{u_0}_{q'}+
\lambda\left\|u\right\|_E^k t^{-\frac{\alpha d (k-1)}{\beta q}+\alpha-\frac{\alpha d k}{\beta}\left(\frac{1}{q'}-\frac{1}{q}\right)}\int_{0}^{1} (1-\tau)^{-\frac{\alpha d (k-1)}{\beta q}+\alpha-1}\tau^{-\frac{\alpha d k}{\beta}\left(\frac{1}{q'}-\frac{1}{q}\right)}d\tau.
\end{align*}
Definition of $q'$ and conditions on $q$, show that
\begin{align*}
\norm{\mathcal{G}u(t)}_q&\lesssim t^{-\frac{\alpha d}{\beta}\left(\frac{1}{q'}-\frac{1}{q}\right)}\norm{u_0}_{q'}+\lambda\left\|u\right\|_E^k t^{\alpha+\frac{\alpha d}{\beta q}-\frac{\alpha d k}{\beta q'}}\frac{\Gamma\left(\alpha-\frac{\alpha d(k-1)}{\beta q}\right)\Gamma\left(1-\frac{\alpha d k}{\beta}\left(\frac{1}{q'}-\frac{1}{q}\right)\right)}{\Gamma\left(1+\alpha+\frac{\alpha d}{\beta q}-\frac{\alpha d k}{\beta q'}\right)}\\
&=C_1t^{-\frac{\alpha d}{\beta}\left(\frac{1}{q'}-\frac{1}{q}\right)}\norm{u_0}_{q'}+C_2\lambda\left\|u\right\|_E^k t^{-\frac{\alpha d}{\beta}\left(\frac{1}{q'}-\frac{1}{q}\right)}\frac{\Gamma\left(\alpha-\frac{\alpha d(k-1)}{\beta q}\right)\Gamma\left(1-\frac{\alpha d k}{\beta}\left(\frac{1}{q'}-\frac{1}{q}\right)\right)}{\Gamma\left(1-\frac{\alpha d}{\beta}\left(\frac{1}{q'}-\frac{1}{q}\right)\right)}
\end{align*}
This proves that $\norm{\mathcal{G}u(t)}_q$ exists for all $t>0$. Furthermore, by defining $\tilde{C_2}:=C_2\frac{\Gamma\left(\alpha-\frac{\alpha d(k-1)}{\beta q}\right)\Gamma\left(1-\frac{\alpha d k}{\beta}\left(\frac{1}{q'}-\frac{1}{q}\right)\right)}{\Gamma\left(1-\frac{\alpha d}{\beta}\left(\frac{1}{q'}-\frac{1}{q}\right)\right)}$, we see that
\begin{align*}
t^{\frac{\alpha d}{\beta}\left(\frac{1}{q'}-\frac{1}{q}\right)}\norm{\mathcal{G}u(t)}_q&\leq C_1\norm{u_0}_{q'}+ \tilde{C_2}\lambda\left\|u\right\|_E^k,
\end{align*}
which consequently establishes $\left\|\mathcal{G}u\right\|_E<\infty$.

Let $D=\{u\in E: \left\|u\right\|_E\leq 2 C_1\norm{u_0}_{q'}\text{ and }u\geq 0\}$. We see that $D$ is a closed set in $E$. Our objective is to demonstrate that $\mathcal{G}$ has a unique fixed point in $D$ by utilizing the contraction mapping principle. This fixed point will act as a super-solution to \eqref{FE:0}. As a result, a solution to \eqref{FE:0} follows directly from Lemma \ref{LemaSuperSol}.

The preceding estimate indicates that if $u\in D$ then
\begin{align*}
\left\|\mathcal{G}u\right\|_E\leq C_1\norm{u_0}_{q'}+\tilde{C_2} \lambda 2^k C_1^k\norm{u_0}_{q'}^k \leq \left( C_1+ \tilde{C_2} \lambda 2^kC_1^k \norm{u_0}_{q'}^{k-1}\right)\norm{u_0}_{q'}
\end{align*}
and thus $\mathcal{G}u \in D$ for small enough $\norm{u_0}_{q'}$. For $u,v\in D$, using the property
\begin{equation*}
|\,|a|^c a-|b|^c b\,|\lesssim |\,a-b\,|(|a|^c+|b|^c)\lesssim |\,a-b\,|(\,|a|+|b|\,)^c,\quad a, b\in\RR,~c>0,
\end{equation*}
and Hölder's inequality, we see that
\begin{align*}
\norm{\mathcal{G}u(t)-\mathcal{G}v(t)}_{q}&\leq C_2\lambda\left\|u+v\right\|_E^{k-1}\left\|u-v\right\|_E\int_{0}^{t} (t-s)^{-\frac{\alpha d(k-1)}{\beta q}+\alpha-1}s^{-\frac{\alpha d k}{\beta}\left(\frac{1}{q'}-\frac{1}{q}\right)}ds\\
&=\tilde{C_2}\lambda\left\|u+v\right\|_E^{k-1}\left\|u-v\right\|_E t^{-\frac{\alpha d}{\beta}\left(\frac{1}{q'}-\frac{1}{q}\right)}
\end{align*}
which implies that
\[
\left\|\mathcal{G}u-\mathcal{G}v\right\|_E\leq C_3 \tilde{C_2} \lambda\norm{u_0}_{q'}^{k-1}\left\|u-v\right\|_E.
\]
As a result, $\mathcal{G}$ acts as a contraction whenever $\norm{u_0}_{q'}$ is sufficiently small. Consequently, there exists a unique fixed point $u\in D$ of $\mathcal{G}$, such that $\mathcal{G} u = u$.

Finally, for some constant $C>0$, we note that
\[
f(s)\leq C \, s^k, \quad s\geq 0,
\]
and therefore, by taking $\lambda=C^{-1}$,
\[
u=\mathcal{G}u\geq \mathcal{F}u\quad a.e.~(t,x)\in (0,T)\times\RR^d
\]
for any $T$ and $\mathcal{F}$ given by \eqref{F:u:S:R:f}.
\end{proof}

\bibliographystyle{plain}
\bibliography{SoVe24bV3Rev}

\begin{thebibliography}{10}

\bibitem{AE87}
J.~Aguirre and M.~Escobedo.
\newblock A {C}auchy problem for $u_t - \triangle u = u^p$ with $0 < p < 1$
  asymptotic behaviour of solutions.
\newblock {\em Annales de la faculté des sciences de Toulouse $5^e$ série},
  8(2):175 -- 203, 1986-1987.

\bibitem{ADK23}
T.~Alinei-Poiana, EH. Dulf, and L.~Kovacs.
\newblock Fractional calculus in mathematical oncology.
\newblock {\em Scientific Reports}, 3(1):10083, 2023.

\bibitem{Aron67}
D.~G. Aronson.
\newblock {Bounds for the fundamental solution of a parabolic equation}.
\newblock {\em Bulletin of the American Mathematical Society}, 73(6):890 --
  896, 1967.

\bibitem{AwMe20}
E.~Awad and R.~Metzler.
\newblock Crossover dynamics from superdiffusion to subdiffusion: Models and
  solutions.
\newblock {\em Fract. Calc. Appl. Anal.}, 23(1):55--102, 2020.

\bibitem{Bel21}
A.~Belmiloudi.
\newblock Cardiac memory phenomenon, time-fractional order nonlinear system and
  bidomain-torso type model in electrocardiology.
\newblock {\em AIMS Mathematics}, 6(1):821--867, 2021.

\bibitem{ClNo79}
Ph. Cl{\'e}ment and J.A. Nohel.
\newblock Abstract linear and nonlinear {V}olterra equations preserving
  positivity.
\newblock {\em SIAM J. Math. Anal.}, 10:365--388, 1979.

\bibitem{EiKo04}
S.-D. Eidelman and A.-N. Kochubei.
\newblock Cauchy problem for fractional diffusion equations.
\newblock {\em J. Differential Equations}, 199(2):211--255, 2004.

\bibitem{JoKo19}
I.~Johnston and V.~Kolokoltsov.
\newblock Green’s {F}unction {E}stimates for {T}ime-{F}ractional {E}volution
  {E}quations.
\newblock {\em Fractal and Fractional}, 3(2), 2019.

\bibitem{KeSVZ16}
J.~Kemppainen, J.~Siljander, V.~Vergara, and R.~Zacher.
\newblock Decay estimates for time-fractional and other non-local in time
  subdiffusion equations in {$\Bbb{R}^d$}.
\newblock {\em Math. Ann.}, 366(3-4):941--979, 2016.

\bibitem{KeSZ17}
J.~Kemppainen, J.~Siljander, and R.~Zacher.
\newblock Representation of solutions and large-time behavior for fully
  nonlocal diffusion equations.
\newblock {\em J. Differential Equations}, 263(1):149--201, 2017.

\bibitem{KiSa04}
A.-A. Kilbas and M.~Saigo.
\newblock {\em {$H$}-transforms}, volume~9 of {\em Analytical Methods and
  Special Functions}.
\newblock Chapman \& Hall/CRC, Boca Raton, FL, 2004.
\newblock Theory and applications.

\bibitem{KiST06}
A.-A. Kilbas, H.-M. Srivastava, and J.-J. Trujillo.
\newblock {\em Theory and applications of fractional differential equations},
  volume 204 of {\em North-Holland Mathematics Studies}.
\newblock Elsevier Science B.V., Amsterdam, 2006.

\bibitem{KiLi16}
K.-H. Kim and S.~Lim.
\newblock Asymptotic behaviors of fundamental solution and its derivatives to
  fractional diffusion-wave equations.
\newblock {\em J. Korean Math. Soc.}, 53(4):929--967, 2016.

\bibitem{Koch90}
A.~N. Kochubei.
\newblock Diffusion of fractional order.
\newblock {\em Differentsial'nye Uravneniya}, 26(4):660--670, 733--734, 1990.

\bibitem{Kolo09}
V.~Kolokoltsov.
\newblock The {L}{\'e}vy--{K}hintchine type operators with variable {L}ipschitz
  continuous coefficients generate linear or nonlinear {M}arkov processes and
  semigroups.
\newblock {\em Probability Theory and Related Fields}, 151(1):95--123, 2011.

\bibitem{Kolo19}
V.~Kolokoltsov.
\newblock {\em Differential equations on measures and functional spaces}.
\newblock Birkh{\"a}user Adv. Texts, Basler Lehrb{\"u}ch. Cham: Birkh{\"a}user,
  2019.

\bibitem{LaRoSi13}
R.~Laister, J.C. Robinson, and M.~Sierzega.
\newblock Non-existence of local solutions for semilinear heat equations of
  {O}sgood type.
\newblock {\em Journal of Differential Equations}, 255(10):3020--3028, 2013.

\bibitem{LW21}
Z.~Lin and H.~Wang.
\newblock Modeling and application of fractional-order economic growth model
  with time delay.
\newblock {\em Fractal and Fractional}, 5(3):74, 2021.

\bibitem{PoVe18}
J.-C. Pozo and V.~Vergara.
\newblock Fundamental solutions and decay of fully non-local problems.
\newblock {\em Discrete Contin. Dyn. Syst.}, 39(1):639--666, 2019.

\bibitem{Pr93}
J.~Pr{\"u}ss.
\newblock {\em Evolutionary {I}ntegral {E}quations and {A}pplications}.
\newblock Birkh{\"a}user Verlag, Basel-Boston-Berlin, 1993.

\bibitem{QS07}
P.~Quittner and P.~Souplet.
\newblock {\em Superlinear {P}arabolic {P}roblems. {B}low-up, {G}lobal
  {E}xistence and {S}teady {S}tates}.
\newblock Birkhäuser Advanced Texts, 2007.

\bibitem{RoSi13}
J.~Robinson and M.~Sierzega.
\newblock Supersolutions for a class of semilinear heat equations.
\newblock {\em Rev Mat Complut}, 26:342--360, 2013.

\bibitem{ScWy89}
W.~R. Schneider and W.~Wyss.
\newblock Fractional diffusion and wave equations.
\newblock {\em J. Math. Phys.}, 30(1):134--144, 1989.

\bibitem{SoVe22}
S.~Solís and V.~Vergara.
\newblock A non-linear stable non-{G}aussian process in fractional time.
\newblock {\em Topological Methods in Nonlinear Analysis}, 59(2B):987–1028,
  Jun. 2022.

\bibitem{SoVe23}
S.~Solís and V.~Vergara.
\newblock Blow-up for a non-linear stable non-{G}aussian process in fractional
  time.
\newblock {\em Fract. Calc. Appl. Anal.}, 2023.

\bibitem{SaTo19}
Sandev T. and Tomovski Z.
\newblock {\em Fractional {E}quations and {M}odels. {T}heory and
  {A}pplications}.
\newblock Springer Nature, 2019.

\bibitem{YGD07}
H.~Ye, J.~Gao, and Y.~Ding.
\newblock A generalized {G}ronwall inequality and its application to a
  fractional differential equation.
\newblock {\em J. Math. Anal. Appl.}, 328:1075 -- 1081, 2007.

\bibitem{ZS15}
Q.~Zhang and H.~Sun.
\newblock The blow-up and global existence of solutions of {C}auchy problems
  for a time fractional diffusion equation.
\newblock {\em Topological Methods in Nonlinear Analysis}, 46:69--92, 2015.

\end{thebibliography}

\end{document}